\documentclass[conference,letterpaper]{IEEEtran}
\usepackage{times}

\usepackage{cite}
\usepackage{graphicx}
\usepackage[cmex10]{amsmath}
\usepackage{array}
\usepackage{subfigure}

\newcommand{\HH}{{\mathcal H}}
\newcommand{\OO}{{\mathcal O}}

\newcommand{\id}{i_d}

\newcommand{\W}{{\Omega}}

\newcommand{\abs}[1]{\left\lvert#1\right\rvert}

\begin{document}

\title{Estimation of Saturation of Permanent-Magnet Synchronous Motors Through an Energy-Based Model}

\author{\IEEEauthorblockN{Al Kassem Jebai\IEEEauthorrefmark{1}, Fran\c{c}ois Malrait\IEEEauthorrefmark{2}, Philippe Martin\IEEEauthorrefmark{1} and Pierre Rouchon\IEEEauthorrefmark{1}}\\
\IEEEauthorblockA{\IEEEauthorrefmark{1}Mines ParisTech, Centre Automatique et Syst\`{e}mes, 60 Bd Saint-Michel, 75272 Paris cedex 06, France\\
Email: \{al-kassem.jebai, philippe.martin, pierre.rouchon\}@mines-paristech.fr}
\IEEEauthorblockA{\IEEEauthorrefmark{2}Schneider Electric, STIE, 33, rue Andr\'{e} Blanchet, 27120 Pacy-sur-Eure, France\\
Email: francois.malrait@schneider-electric.com}}

\maketitle
\begin{abstract}
We propose a parametric model of the saturated Permanent-Magnet Synchronous Motor (PMSM) together with an estimation method of the magnetic parameters. The model is based on an energy function which simply encompasses the saturation effects. Injection of fast-varying pulsating voltages and measurements of the resulting current ripples then permit to identify the magnetic parameters by linear least squares. Experimental results on a surface-mounted PMSM and an interoir magnet PMSM illustrate the relevance of the approach.
\end{abstract}

\section{Introduction}

Sensorless control of Permanent-Magnet Synchronous Motors (PMSM) at low velocity remains a challenging task. Most of the existing control algorithms rely on the motor saliency, both geometric and saturation-induced, for extracting the rotor position from the current measurements through high-frequency signal injection~\cite{Ji_Sul_2003,CorleL1998IToIA}. However some magnetic saturation effects such as cross-coupling and permanent magnet demagnetization can introduce large errors on the rotor position estimation~\cite{GugliPV2006IToIA,Bi_bo_2009}. These errors decrease the performance of the controller. In some cases they may cancel the rotor total saliency and lead to instability. It is thus important to correctly model the magnetic saturation effects, which is usually done through d-q magnetizing curves (flux versus current). These curves are usually found either by finite element analysis FEA or experimentally by integration of the voltage equation~\cite{StumSDHT2003IToIA, Stumb2005}. This provides a good way to characterize the saturation effects and can be used to improve the sensorless control of the PMSM~\cite{Zh_li_2007,Gugl_2004}. However the FEA or the integration of the voltage equation methods are not so easy to implement and do not provide an explicit model of the saturated PMSM.

In this paper a simple parametric model of the saturated PMSM is introduced (section~\ref{sec:model}); it is based on an energy function~\cite{BasicMR2010IToAC,Dur_2011} which simply encompasses the saturation and cross-magnetization effects.
In section~\ref{sec:Identi} a simple estimation method of the magnetic parameters is proposed and rigorously justified: fast-varying pulsating voltages are impressed to the motor with rotor locked; they create current ripples from which the magnetic parameters are estimated by linear least squares. In section~\ref{sec:experiment} experimental results on two kinds of motors (with surface-mounted and interior magnets) illustrate the relevance of the approach.

\section{An energy-based model for the saturated PMSM}\label{sec:model}
\subsection{Energy-based model}
The electrical subsystem of a two-axis PMSM expressed in the synchronous $d-q$ frame reads
\begin{align}
\label{eq:elecd}\frac{d\phi_d}{dt} &=u_d-Ri_d+\frac{d\theta}{dt}\phi_q\\
\label{eq:elecq}\frac{d\phi_q}{dt} &=u_q-Ri_q-\frac{d\theta}{dt}(\phi_d+\phi_m),
\end{align}
where $\phi_d,\phi_m$ are the direct-axis flux linkages due to the current excitation and to the permanent magnet, and $\phi_q$ is the quadrature-axis flux linkage; $u_d,u_q$ are the impressed voltages and $i_d,i_q$ are the currents; $\theta$ is the rotor (electrical) position and $R$ is the stator resistance.
The currents can be expressed in function of the flux linkages thanks to a suitable energy function~$\HH(\phi_d,\phi_q)$ by
\begin{align}
   \label{eq:CurrentFluxd}i_d &=\partial_1\HH(\phi_d,\phi_q)\\
   \label{eq:CurrentFluxq}i_q &=\partial_2\HH(\phi_d,\phi_q),
\end{align}
where $\partial_k\HH$ denotes the partial derivative w.r.t. the $k^{th}$ variable, see~\cite{BasicMR2010IToAC,Dur_2011}; without loss of generality $\HH(0,0)=0$.

For an unsaturated PMSM this energy function reads
\begin{align*}
\HH_l(\phi_d,\phi_q) &=\frac{1}{2L_d}\phi^2_d+\frac{1}{2L_q}\phi^2_q
\end{align*}
where $L_d$ and $L_q$ are the motor self-inductances, and we recover the usual linear relations
\begin{align*}
    i_d &=\partial_1\HH(\phi_d,\phi_q) =\frac{\phi_d}{L_d}\\
    i_q &=\partial_2\HH(\phi_d,\phi_q) =\frac{\phi_q}{L_q}.
\end{align*}

Notice the expression for $\HH$ should respect the symmetry of the PMSM w.r.t the direct axis, i.e.
\begin{equation}\label{eq:sym}
    \HH(\phi_d,-\phi_q)=\HH(\phi_d,\phi_q),
\end{equation}
which is obviously the case for~$\HH_l$.
Indeed, \eqref{eq:elecd}-\eqref{eq:elecq} is left unchanged by the transformation
$$(\phi_d',u_d',i_d',\phi_q',u_q',i_q',\theta'):=(\phi_d,u_d,i_d,-\phi_q,-u_q,-i_q,-\theta);$$
this implies
\begin{align*}
    \partial_1\HH(\phi_d',\phi_q') &=\partial_1\HH(\phi_d,\phi_q)\\
    \partial_2\HH(\phi_d',\phi_q') &=-\partial_2\HH(\phi_d,\phi_q),
\end{align*}
i.e.
\begin{align*}
    \partial_1\HH(\phi_d,-\phi_q) &=\partial_1\HH(\phi_d,\phi_q)\\
    \partial_2\HH(\phi_d,-\phi_q) &=-\partial_2\HH(\phi_d,\phi_q).
\end{align*}
Therefore
\begin{align*}
    \frac{d\HH}{d\phi_d}(\phi_d,-\phi_q) &=\partial_1\HH(\phi_d,-\phi_q)\\
    &=\partial_1\HH(\phi_d,\phi_q)\\
    &=\frac{d\HH}{d\phi_d}(\phi_d,\phi_q)\\
    \frac{d\HH}{d\phi_q}(\phi_d,-\phi_q) &=-\partial_2\HH(\phi_d,-\phi_q)\\
    &=\partial_2\HH(\phi_d,\phi_q)\\
    &=\frac{d\HH}{d\phi_q}(\phi_d,\phi_q).
\end{align*}
Integrating these relations yields
\begin{align*}
    \HH(\phi_d,-\phi_q) &=\HH(\phi_d,\phi_q)+c_d(\phi_q)\\
    \HH(\phi_d,-\phi_q) &=\HH(\phi_d,\phi_q)+c_q(\phi_d),
\end{align*}
where $c_d,c_q$ are functions of only one variable. But this makes sense only if $c_d(\phi_q)=c_q(\phi_d)=c$ with $c$ constant. Since $H(0,0)=0$, $c=0$, which yields~\eqref{eq:sym}.

\subsection{Parametric description of magnetic saturation}\label{sec:paramsat}
\begin{figure}[t!]\centering
\subfigure[$\phi_d(i_d,i_q=Constant)$]{\includegraphics[width=\columnwidth]{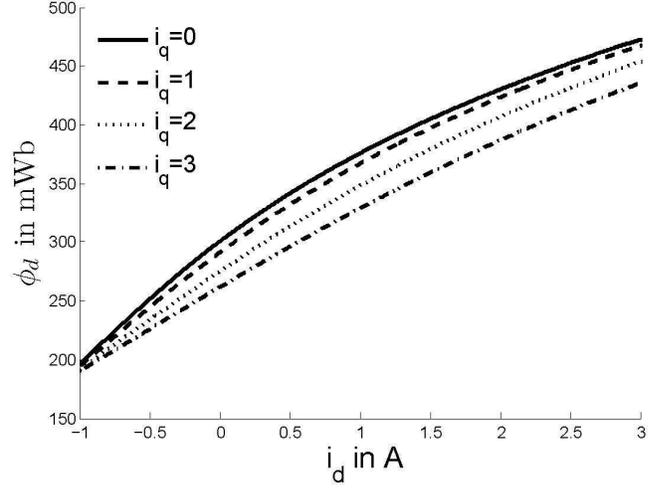}\label{fig:idphid}}
\subfigure[$\phi_q(i_d=Constant,i_q)$]{\includegraphics[width=\columnwidth]{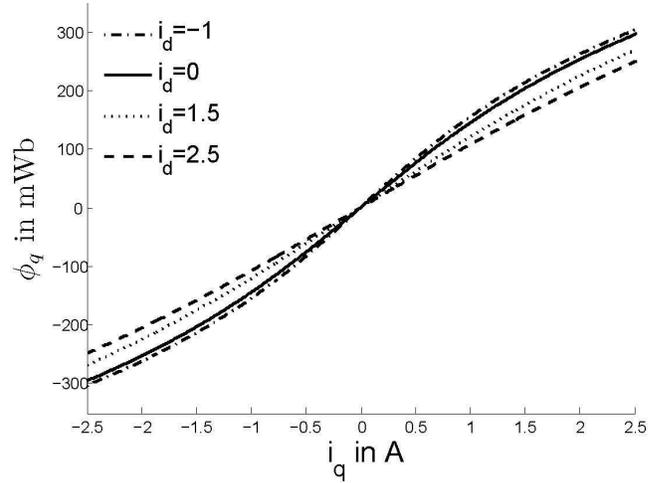}\label{fig:iqphiq}}
\caption{Flux-current magnetization curves (IPM)}\label{fig:iphi}
\end{figure}

Magnetic saturation can be accounted for by considering a more complicated magnetic energy function~$\HH$, having $\HH_l$ for quadratic part but including also higher-order terms. From experiments saturation effects are well captured by considering only third- and fourth-order terms, hence
\begin{multline*}
\HH(\phi_d,\phi_q)=\HH_l(\phi_d,\phi_q)\\
+\sum_{i=0}^3\alpha_{3-i,i}\phi_d^{3-i}\phi_q^i+\sum_{i=0}^4\alpha_{4-i,i}\phi_d^{4-i}\phi_q^i.
\end{multline*}
This is a perturbative model where the higher-order terms appear as corrections of the dominant term~$\HH_l$.
The $9$ coefficients $\alpha_{ij}$
together with $L_d$, $L_q$ are motor dependent. But \eqref{eq:sym} implies $\alpha_{2,1}=\alpha_{0,3}=\alpha_{3,1}=\alpha_{1,3}=0$, so that the energy function eventually reads
\begin{multline}\label{eq:EnerSat}
\HH(\phi_d,\phi_q) =\HH_l(\phi_d,\phi_q) +\alpha_{3,0}\phi_d^3+\alpha_{1,2}\phi_d\phi_q^2 \\
+\alpha_{4,0}\phi_d^4+\alpha_{2,2}\phi_d^2\phi_q^2+\alpha_{0,4}\phi_q^4.
\end{multline}
From~\eqref{eq:CurrentFluxd}-\eqref{eq:CurrentFluxq} and~\eqref{eq:EnerSat} the currents are then explicitly given by
\begin{align}
 \notag i_d &=\partial_1\HH(\phi_d,\phi_q)\\
 \label{eq:id}&=\frac{\phi_d}{L_d}+3\alpha_{3,0}\phi_d^2+\alpha_{1,2}\phi_q^2 +4\alpha_{4,0}\phi_d^3+2\alpha_{2,2}\phi_d\phi_q^2\\
 \notag i_q &=\partial_2\HH(\phi_d,\phi_q)\\
 \label{eq:iq}&=\frac{\phi_q}{L_q}+2\alpha_{1,2}\phi_d\phi_q+2\alpha_{2,2}\phi_d^2\phi_q+4\alpha_{0,4}\phi_q^3,
\end{align}
which are the flux-current magnetization curves. Fig.~\ref{fig:iphi} shows examples of these curves in the more familiar presentation of fluxes w.r.t currents obtained by numerically inverting~\eqref{eq:CurrentFluxd}-\eqref{eq:CurrentFluxq}; the motor is the IPM of section~\ref{sec:experiment}.

The model of the saturated PMSM is thus given by \eqref{eq:elecd}-\eqref{eq:elecq} and~\eqref{eq:id}-\eqref{eq:iq}. It is in state form with $\phi_d,\phi_q$ as state variables. The magnetic saturation effects are represented by the $5$ additional parameters $\alpha_{3,0}, \alpha_{1,2}, \alpha_{4,0}, \alpha_{2,2}, \alpha_{0,4}$.

\subsection{Model with $i_d,i_q$ as state variables}
The model of the saturated PMSM is often expressed with $i_d,i_q$ as state variables, e.g.~\cite{StumSDHT2003IToIA}. Starting with flux-current magnetization curves in the form
\begin{align}
    \label{eq:FluxCurrentd}\phi_d &=\Phi_d(i_d,i_q)\\
    \label{eq:FluxCurrentq}\phi_q &=\Phi_q(i_d,i_q)
\end{align}
and differentiating w.r.t time,
\eqref{eq:elecd}-\eqref{eq:elecq} then becomes
\begin{align*}
L_{dd}(i_d,i_q)\frac{di_d}{dt}+L_{dq}(i_d,i_q)\frac{di_q}{dt} &=u_d-Ri_d+\frac{d\theta}{dt}\phi_q\\
L_{qd}(i_d,i_q)\frac{di_d}{dt}+L_{qq}(i_d,i_q)\frac{di_q}{dt} &=u_q-Ri_q-\frac{d\theta}{dt}(\phi_d+\phi_m),
\end{align*}
where
\begin{align*}
   \begin{pmatrix}L_{dd}(i_d,i_q) & L_{dq}(i_d,i_q)\\ L_{qd}(i_d,i_q)& L_{qq}(i_d,i_q)\end{pmatrix}
   =&\begin{pmatrix}\partial_1\Phi_d(i_d,i_q) & \partial_2\Phi_d(i_d,i_q)\\
   \partial_1\Phi_q(i_d,i_q) & \partial_2\Phi_q(i_d,i_q)\end{pmatrix}.
\end{align*}

Though not always acknowledged $L_{dq}$ and $L_{qd}$ should be equal. Indeed, plugging \eqref{eq:CurrentFluxd}-\eqref{eq:CurrentFluxq} into \eqref{eq:FluxCurrentd}-\eqref{eq:FluxCurrentq} gives
\begin{align*}
   \phi_d &=\Phi_d\bigl(\partial_1\HH(\phi_d,\phi_q),\partial_2\HH(\phi_d,\phi_q)\bigr)\\
   \phi_q &=\Phi_q\bigl(\partial_1\HH(\phi_d,\phi_q),\partial_2\HH(\phi_d,\phi_q)\bigr).
\end{align*}
Taking the total derivative of both sides of these equations w.r.t. $\phi_d$ and $\phi_q$ then yields
\begin{align*}
   \begin{pmatrix}1&0\\0&1\end{pmatrix}
   =&\begin{pmatrix}L_{dd}\partial_{11}\HH+L_{dq}\partial_{12}\HH & L_{dd}\partial_{21}\HH+L_{dq}\partial_{22}\HH\\
   L_{qd}\partial_{11}\HH+L_{qq}\partial_{12}\HH & L_{qd}\partial_{21}\HH+L_{qq}\partial_{22}\HH\end{pmatrix}\\
   =&\begin{pmatrix}L_{dd} & L_{dq}\\ L_{qd} & L_{qq}\end{pmatrix}
   \begin{pmatrix}\partial_{11}\HH & \partial_{21}\HH\\ \partial_{12}\HH & \partial_{22}\HH\end{pmatrix}.
\end{align*}
Since $\partial_{12}\HH=\partial_{21}\HH$ the second matrix in the last line is symmetric, hence the first; in other words $L_{dq}=L_{qd}$.

To do that with the model of section~\ref{sec:paramsat} the nonlinear equations \eqref{eq:id}-\eqref{eq:iq} must be inverted. Rather than doing that exactly, we take advantage of the fact the coefficients $\alpha_{i,j}$ are experimentally small. At first order w.r.t. the~$\alpha_{i,j}$ we obviously have $\phi_d=L_di_d+\OO(\abs{\alpha_{i,j}})$ and $\phi_q=L_qi_q+\OO(\abs{\alpha_{i,j}})$. Plugging these expressions into~\eqref{eq:id}-\eqref{eq:iq} we easily find
\begin{align}
 \notag\phi_d &=L_d\bigl(i_d-3\alpha_{3,0}L_d^2i_d^2 -\alpha_{1,2}L_q^2i_q^2-4\alpha_{4,0}L_d^3i_d^3\\
 \label{eq:phid}&\quad-2\alpha_{2,2}L_dL_q^2i_di_q^2\bigr)+\OO(\abs{\alpha_{i,j}}^2)\\
 \notag\phi_q &=L_q\bigl(i_q-2\alpha_{1,2}L_dL_qi_di_q-2\alpha_{2,2}L_d^2L_qi_d^2i_q\\
 \label{eq:phiq}&\quad-4\alpha_{0,4}L_q^3i_q^3\bigr)+\OO(\abs{\alpha_{i,j}}^2).
\end{align}
Finally,
\begin{align*}
   L_{dd}(i_d,i_q) &=L_d\bigl(1-6\alpha_{3,0}L_di_d-12\alpha_{4,0}L_d^2i_d^2-2\alpha_{2,2}L_q^2i_q^2\bigr)\\
   L_{dq}(i_d,i_q) &=L_{qd}(i_d,i_q)=-2L_dL_q^2i_q(\alpha_{1,2}+2\alpha_{2,2}L_di_d)\\
   L_{qq}(i_d,i_q) &=L_q\bigl(1-2\alpha_{1,2}L_di_d-2\alpha_{2,2}L_d^2i_d^2-12\alpha_{0,4}L_q^2i_q^2\bigr).
\end{align*}

\section{A procedure for estimating the magnetic parameters}\label{sec:Identi}
\subsection{Principle}\label{sec:principle}
\begin{figure}[t!]
\centerline{\includegraphics[width=1\columnwidth]{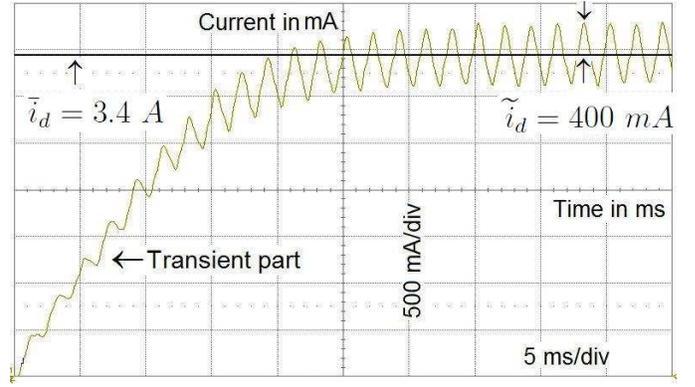}}
\caption{Experimental illustration of equation~\eqref{eq:iprocd}: time response of~$i_d$}\label{fig:HFcurrent}
\end{figure}
To estimate the $7$ magnetic parameters in the model, we propose a procedure which is rather easy to implement and reliable. With the rotor locked in the position~$\theta=0$, we inject fast-varying pulsating voltages
\begin{align}
   \label{eq:uprocd}u_d(t) &=\bar u_d+\widetilde u_df(\W t)\\
   \label{eq:uprocq}u_q(t) &=\bar u_q+\widetilde u_qf(\W t),
\end{align}
where $\bar u_d,\bar u_q,\widetilde u_d,\widetilde u_q,\W$ are constant and $f$ is a periodic function with zero mean. The pulsation $\W$ is chosen large enough w.r.t. the motor electric time constant. It can then be shown, see section~\ref{sec:justification}, that after an initial transient
\begin{align}
    \label{eq:iprocd}i_d(t) &=\bar i_d+\widetilde i_dF(\W t)+\OO(\tfrac{1}{\W^2})\\
    \label{eq:iprocq}i_q(t) &=\bar i_q+\widetilde i_qF(\W t)+\OO(\tfrac{1}{\W^2}),
\end{align}
where $\bar i_d=\frac{\bar u_d}{R},\bar i_q=\frac{\bar u_q}{R},\widetilde i_d,\widetilde i_q$ are constant and $F$ is the primitive of~$f$ with zero mean ($F$ has clearly the same period as~$f$); fig.~\ref{fig:HFcurrent} shows for instance the current $i_d$ obtained for the SPM of section~\ref{sec:experiment} when starting from $\id(0)=0$ and applying a square signal $u_d$ with $\Omega=500Hz$, $\bar u_d=23V$ and~$\widetilde u_d=30V$.
On the other hand using the saturation model the amplitudes $\widetilde i_d,\widetilde i_q$ of the current ripples turn out to be
\begin{align}
    \notag\widetilde i_d &=\frac{1}{\W}\Bigl(\frac{\widetilde u_d}{L_d}
    +2\alpha_{2,2}L_q\bar i_q(2L_d\bar i_d\widetilde u_q+L_q\bar i_q\widetilde u_d)\\
    \label{eq:tid}&\quad+12\alpha_{4,0}L_d^2\bar i_d^2\widetilde u_d
    +6\alpha_{3,0}L_d\bar i_d\widetilde u_d+2\alpha_{1,2}L_q\bar i_q\widetilde u_q\Bigr)\\
    \notag\widetilde i_q &=\frac{1}{\W}\Bigl(\frac{\widetilde u_q}{L_q}
    +2\alpha_{2,2}L_d\bar i_d(2L_q\bar i_q\widetilde u_d+L_d\bar i_d\widetilde u_q)\\
    \label{eq:tiq}&\quad+12\alpha_{0,4}L_q^2\bar i_q^2\widetilde u_q
    +2\alpha_{1,2}(L_d\bar i_d\widetilde u_q+L_q\bar i_q\widetilde u_d)\Bigr).
\end{align}
As $\widetilde i_d,\widetilde i_q$ can easily be measured experimentally, these expressions provide a means to identify the magnetic parameters from experimental data obtained with various values of~$\bar u_d,\bar u_q,\widetilde u_d,\widetilde u_q$.

\subsection{Estimation of the parameters}
Since combinations of the magnetic parameters always enter~\eqref{eq:tid}-\eqref{eq:tiq} linearly, they can be estimated by simple linear least squares; moreover by suitably choosing $\bar u_d,\bar u_q,\widetilde u_d,\widetilde u_q$, the whole least squares problem for the $7$ parameters can be split into several problems involving fewer parameters:
\begin{itemize}
 \item with $\bar u_d=\bar u_q=0$, hence $\bar i_d=\bar i_q=0$, and $\widetilde u_d=0$ (resp.~$\widetilde u_q=0$) equation~\eqref{eq:tid} (resp. equation~\eqref{eq:tiq}) reads
\begin{equation}\label{eq:LdLq}
    L_d=\frac{1}{\W}\frac{\widetilde u_d}{\widetilde i_d}
    \qquad\Bigl(\text{resp.~~}L_q=\frac{1}{\W}\frac{\widetilde u_q}{\widetilde i_q}\Bigr)
\end{equation}
 \item with $\bar u_q=0$, hence $\bar i_q=0$, and $\widetilde u_q=0$, \eqref{eq:tid} reads
\begin{equation}\label{eq:saturation1}
    \widetilde i_d=\frac{\widetilde u_d}{\W}\left(\frac{1}{L_d}+6\alpha_{3,0}L_d\bar i_d+12\alpha_{4,0}L_d^2\bar i_d^2\right).
\end{equation}
Notice \eqref{eq:tiq} reads $\widetilde i_q=0$ hence provides no information
 \item with $\bar u_d=0$, hence $\bar i_d:=0$, and $\widetilde u_q=0$, \eqref{eq:tid}-\eqref{eq:tiq} read
\begin{align}
   \label{eq:saturation3d}\widetilde i_d  &=\frac{\widetilde u_d}{\W}\Bigl(\frac{1}{L_d}+2\alpha_{2,2}L_q^2\bar i_q^2\Bigr)\\
   \label{eq:saturation3q}\widetilde i_q & =\frac{2\widetilde u_d}{\W}\alpha_{1,2}L_q\bar i_q
\end{align}
 \item with $\bar u_d=0$, hence $\bar i_d:=0$, and $\widetilde u_d=0$, \eqref{eq:tid}-\eqref{eq:tiq} read
\begin{align}
  \label{eq:saturation4d}\widetilde i_d &=\frac{2\widetilde u_q}{\W}\alpha_{1,2}L_q\bar i_q\\
  \label{eq:saturation4q}\widetilde i_q &=\frac{\widetilde u_q}{\W}\Bigl(\frac{1}{L_q}+12\alpha_{0,4}L_q^2\bar i_q^2\Bigr).
\end{align}
\end{itemize}

$L_d$ (resp.~$L_q$) is then immediately determined from~\eqref{eq:LdLq}; $\alpha_{3,0}$ and $\alpha_{4,0}$ are jointly estimated by least squares from~\eqref{eq:saturation1}; $\alpha_{2,2}$, $\alpha_{1,2}$ and $\alpha_{0,4}$ are separately estimated by least squares from respectively \eqref{eq:saturation3d}, \eqref{eq:saturation3q}-\eqref{eq:saturation4d} and~\eqref{eq:saturation4q}.

\subsection{Justification of section~\ref{sec:principle}}\label{sec:justification}
The assertions of section~\ref{sec:principle} can be rigorously justified by a straightforward application of second-order averaging of differential equations~\cite[p.~40]{SandeVM2007book}. Indeed the electrical subsystem \eqref{eq:elecd}-\eqref{eq:elecq} with locked rotor (i.e. $\frac{d\theta}{dt}=0$) and input voltages~\eqref{eq:uprocd}-\eqref{eq:uprocd} reads when setting $\tau=\W t$
\begin{align}
    \label{eq:dynaved}\frac{d\phi_d}{d\tau} & =\frac{1}{\W}\bigl(\bar u_d+\widetilde u_df(\tau)-Ri_d(\phi_d,\phi_q)\bigr)\\
    \label{eq:dynaveq}\frac{d\phi_q}{d\tau} & =\frac{1}{\W}\bigl(\bar u_q+\widetilde u_qf(\tau)-Ri_q(\phi_d,\phi_q)\bigr).
\end{align}
This system is in the so-called standard form for averaging, with a right hand-side periodic in~$\tau$ and $\frac{1}{\W}$ as a small parameter. Therefore its solution is given by
\begin{align}
   \label{eq:phiprocd}\phi_d(\tau) &=\phi_d^0(\tau)+\frac{\widetilde u_d}{\W}F(\tau)+\OO(\tfrac{1}{\W^2})\\
   \label{eq:phiprocq}\phi_q(\tau) &=\phi_q^0(\tau)+\frac{\widetilde u_q}{\W}F(\tau) +\OO(\tfrac{1}{\W^2}),
\end{align}
where $(\phi_d^0,\phi_q^0)$ is the solution of the system
\begin{align*}
\frac{d\phi_d^0}{dt} &=\bar u_d-Ri_d(\phi_d^0,\phi_q^0)\\
\frac{d\phi_q^0}{dt} &=\bar u_q-Ri_q(\phi_d^0,\phi_q^0)
\end{align*}
obtained by averaging the right-hand side of~\eqref{eq:dynaved}-\eqref{eq:dynaveq}. After an initial transient
$\bigl(\phi_d^0(\tau),\phi_q^0(\tau)\bigr)$ asymptotically reaches the constant value $(\bar\phi_d,\bar\phi_q)$ determined by $\bar u_d=Ri_d(\bar\phi_d,\bar\phi_q)$ and $\bar u_q=Ri_q(\bar\phi_d,\bar\phi_q)$.

Plugging~\eqref{eq:phiprocd}-\eqref{eq:phiprocq} with $t=\frac{\tau}{\W}$ into~\eqref{eq:id}-\eqref{eq:iq}, and expanding along powers of $\frac{1}{\W}$ then yields
\begin{align*}
i_d(t) &=\bar i_d+\frac{F(\W t)}{\W}\biggr(\frac{\widetilde u_d}{L_d}+6\alpha_{3,0}\bar\phi_d\widetilde u_d+2\alpha_{1,2}\bar\phi_q\widetilde u_q\\
 &\quad+12\alpha_{4,0}\bar\phi_d^2\widetilde u_d+2\alpha_{2,2}(2\bar\phi_d\bar\phi_q\widetilde u_q+\bar\phi_q^2\widetilde u_d)\biggl)+\OO(\tfrac{1}{\W^2})\\
i_q(t) &=\bar i_q+\frac{F(\W t)}{\W}\biggl(\frac{\widetilde u_q}{L_q}+2\alpha_{1,2}(\bar\phi_d\widetilde u_q+\bar\phi_q\widetilde u_d)\\
&\quad+2\alpha_{2,2}(2\bar\phi_d\bar\phi_q\widetilde u_d+\bar\phi_d^2\widetilde u_q)
+12\alpha_{0,4}\bar\phi_q^2\widetilde u_q\biggr)
+\OO(\tfrac{1}{\W^2}),
\end{align*}
where $\bar i_d=i_d(\bar\phi_d,\bar\phi_q)$ and $\bar i_q=i_q(\bar\phi_d,\bar\phi_q)$
There remains to express $\bar\phi_d,\bar\phi_q$ in function of $\bar i_d,\bar i_q$. Rather than exactly inverting the nonlinear equations~\eqref{eq:id}-\eqref{eq:iq}, we take advantage of the fact the coefficients $\alpha_{i,j}$ are experimentally small. At first order w.r.t. the~$\alpha_{i,j}$ we have $\phi_d=L_di_d+\OO(\abs{\alpha_{i,j}})$ and $\phi_q=L_qi_q+\OO(\abs{\alpha_{i,j}})$. Using this in the previous equations and neglecting $\OO(\tfrac{1}{\W^2})$ and $\OO(\abs{\alpha_{i,j}}^2)$ terms we eventually find~\eqref{eq:iprocd}-\eqref{eq:tiq}. Using directly \eqref{eq:phid}-\eqref{eq:phiq} yields of course the same result.
\section{Experimental Results}\label{sec:experiment}
\subsection{Experimental setup}
The methodology of section~\ref{sec:Identi} is tested on an interior magnet PMSM (IPM) and a surface-mounted PMSM (SPM) with rated parameters listed below. 
The setup consists of an industrial inverter with a $400V$ DC bus and a $4kHz$ PWM switching frequency, 3 dSpace boards (DS1005 PPC Board, DS2003 A/D Board, DS4002 Timing and Digital I/O Board) and a host PC. The measurements were sampled also at $4kHz$.
%
\begin{center}
\begin{tabular}{| l | l | l | }
  \firsthline
   & IPM & SPM  \\
  \hline
  Pole pairs & 6 & 2 \\
  \hline
  Rated power & $200~W$ & $1200~W$ \\
  \hline
  Rated current & $1.2~A$ & $3.4~A$ \\
  \hline
  Rated speed & $1800~rpm$ & $400~rpm$ \\
  \hline
  Rated torque & $1.06~N.m$ & $29~N.m$ \\
  \hline
  Resistance & $12.15~\W$ & $6.69~\W$ \\
  \lasthline
\end{tabular}
\end{center}

\subsection{Experimental results}
With the rotor locked in the position~$\theta=0$, a square wave voltage with frequency $\W=500Hz$ and constant amplitude $\widetilde u_d$ or $\widetilde u_q$ ($30V$ for the IPM, $40V$ for the SPM) is applied to the motor. But for the determination of $L_d,L_q$ where $\bar u_d=\bar u_q=0$, several runs are performed with various $\bar u_d$ (resp.~$\bar u_q$) such that $\bar i_d$ (resp.~$\bar i_q$) ranges from $-2A$ to $+2A$ with a $0.3A$ increment (IPM), or from $-8A$ to $8A$ with a $0.5A$ increment (SPM). The estimated parameters are listed below; the uncertainty in the estimation stems from a $\pm10mA$ uncertainty in the current measurements.
\begin{center}
\begin{tabular}{|l | l | l | }
\firsthline
     & IPM & SPM \\
  \hline
  $L_d~(mH$) & $91.9\pm5$ & $155.4\pm10$ \\
   \hline
  $L_q~(mH)$ & $45.8\pm1$ & $58.6\pm2$ \\
   \hline
  $\alpha_{3,0}~(A.Wb^{-2})$ & $7.70\pm0.11$ & $5.01\pm0.11$ \\
   \hline
  $\alpha_{1,2}~(A.Wb^{-2})$ & $5.35\pm0.61$ & $4.83\pm0.27$ \\
   \hline
  $\alpha_{4,0}~(A.Wb^{-3})$ & $19.42\pm1.34$ & $1.83\pm0.28$ \\
   \hline
  $\alpha_{2,2}~(A.Wb^{-3})$ & $22.18\pm2.80$ & $8.76\pm1.03$ \\
   \hline
  $\alpha_{0,4}~(A.Wb^{-3})$ & $6.62\pm0.42$ & $1.18\pm0.17$ \\
  \lasthline
\end{tabular}
\end{center}

\begin{figure}[t!]
\centering
\subfigure[IPM]{\includegraphics[height=67mm]{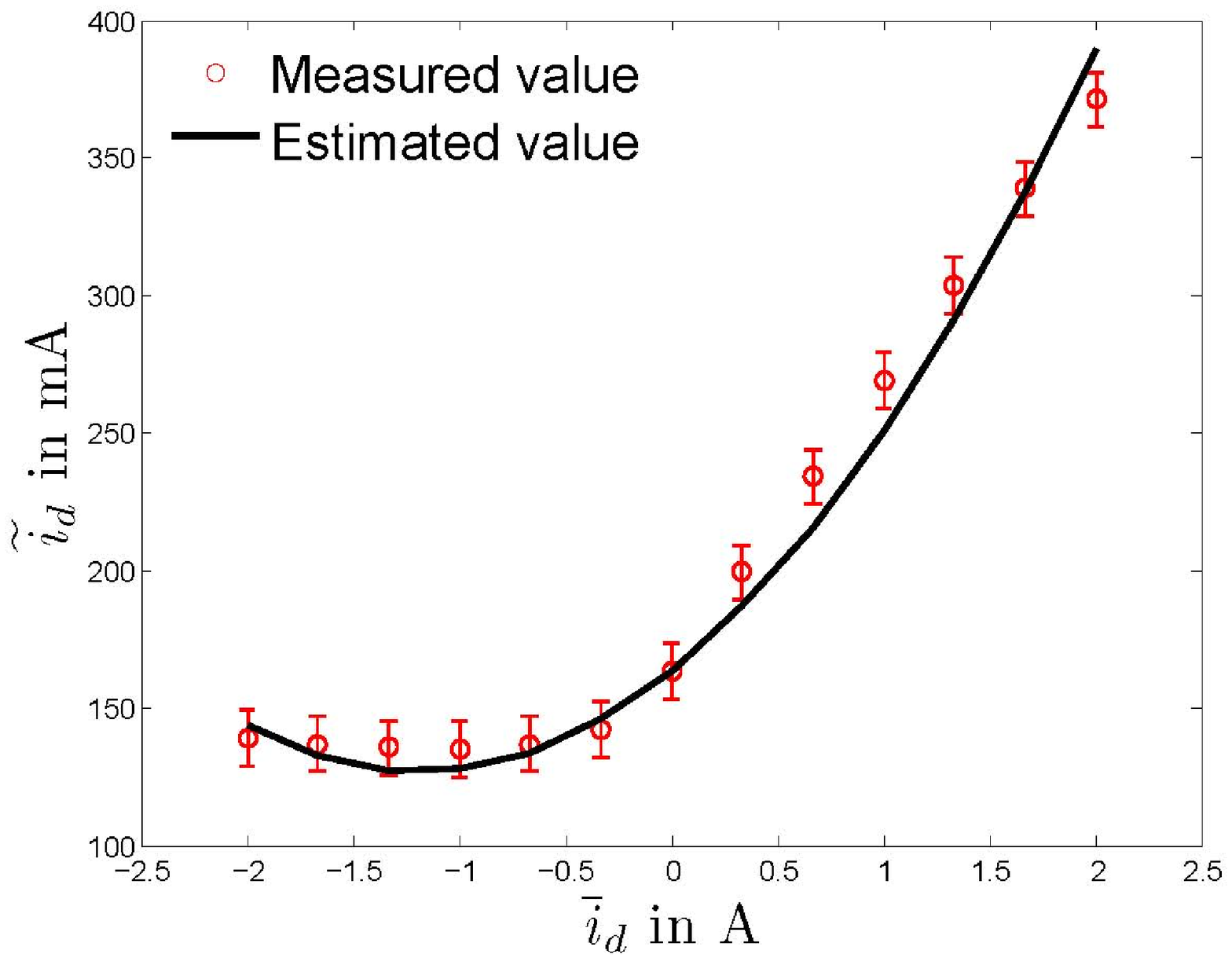}\label{fig:Id_id_ud2}}
\subfigure[SPM]{ \includegraphics[height=67mm]{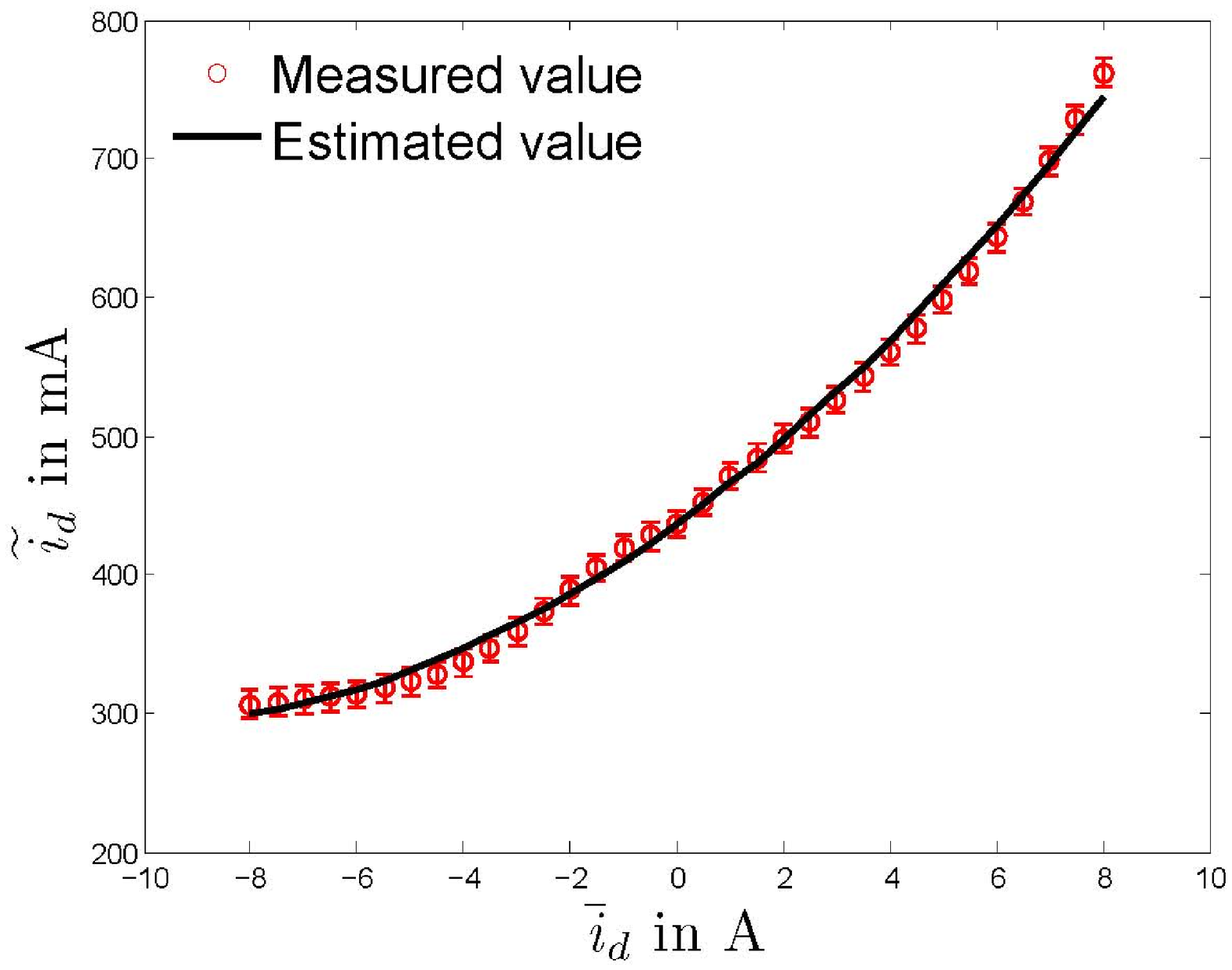}\label{fig:Id_id_ud1}}
\caption{Measured values (circles) and fitted curve (solid line) for~\eqref{eq:saturation1}.}\label{fig:Id_id_ud}
\end{figure}

\begin{figure}[t!]
\centering
\subfigure[IPM]{\includegraphics[height=67mm]{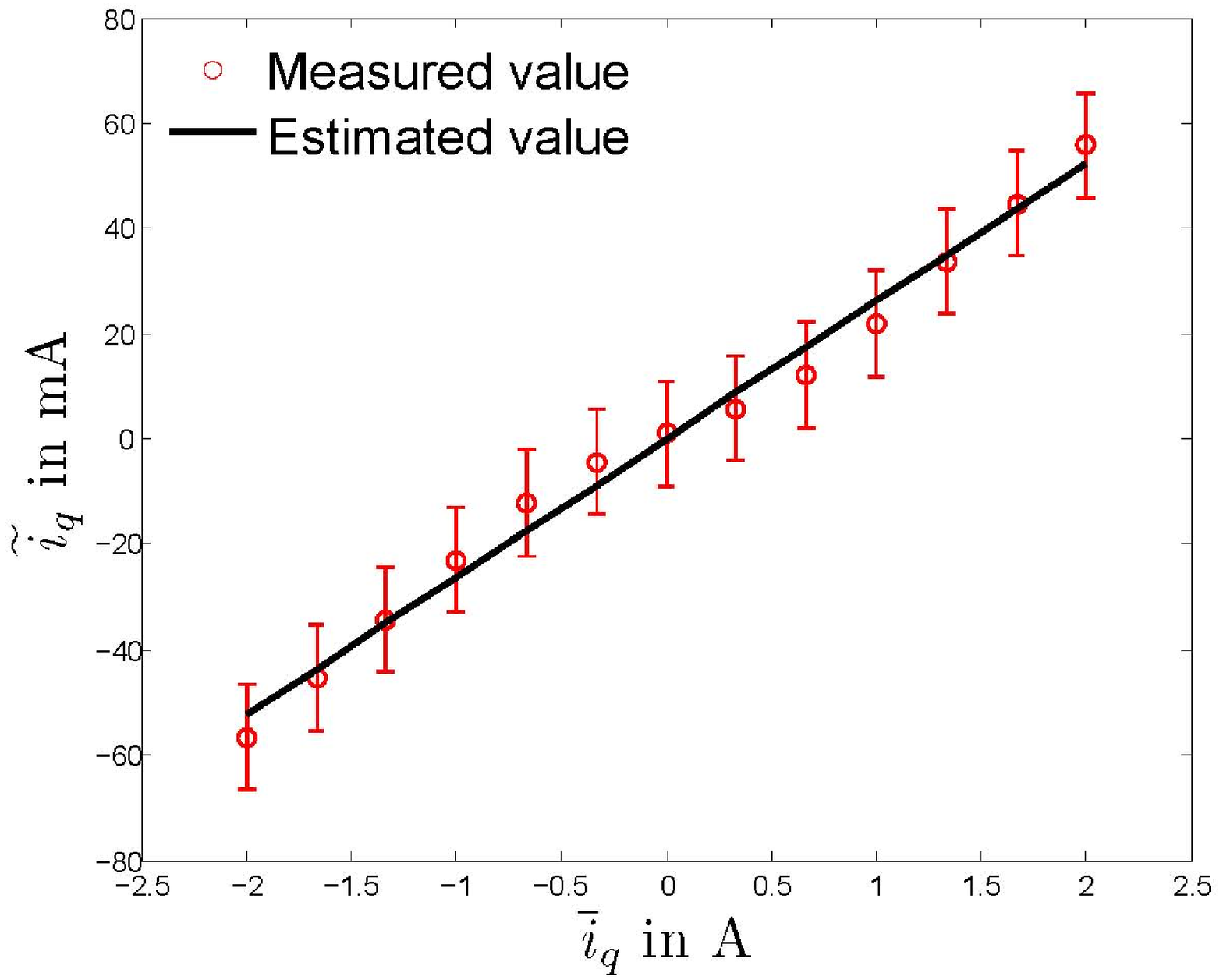}\label{fig:Iq_iq_ud2}}
\subfigure[SPM]{\includegraphics[height=67mm]{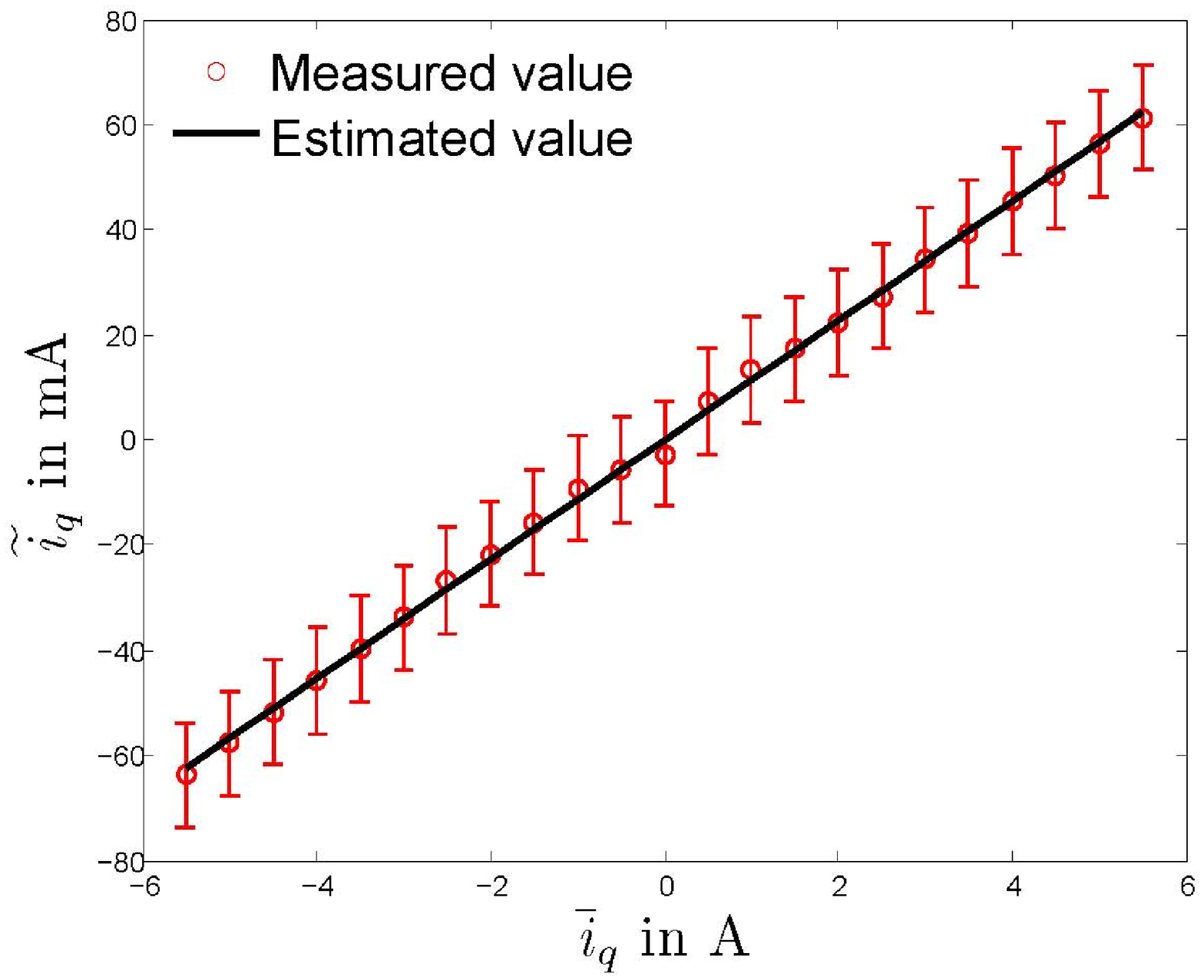}\label{fig:Iq_iq_ud1}}
\caption{Measured values (circles) and fitted curve (solid line) for~\eqref{eq:saturation3q}.}\label{fig:Iq_iq_ud}
\end{figure}

The good agreement between the fitted curves and the measurements is demonstrated for instance for~\eqref{eq:saturation1} on~Fig.~\ref{fig:Id_id_ud} and for~\eqref{eq:saturation3q} on~Fig.~\ref{fig:Iq_iq_ud}. Notice \eqref{eq:saturation1} illustrates saturation on a single axis, while \eqref{eq:saturation3q} illustrates cross-saturation.

%
%

\begin{figure}[t!]
\centering
\subfigure[Interior-magnet PMSM]{
       \includegraphics[width=\columnwidth]{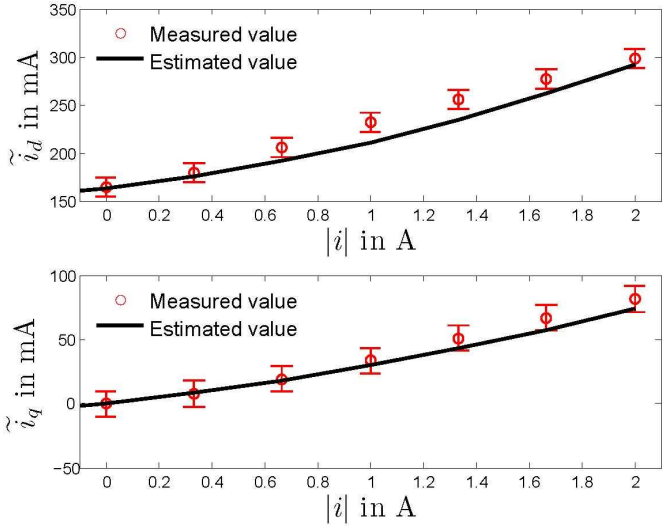}
   \label{fig:Id_i60_ud2}
 }
 \subfigure[Surface mounted PMSM]{
       \includegraphics[width=\columnwidth]{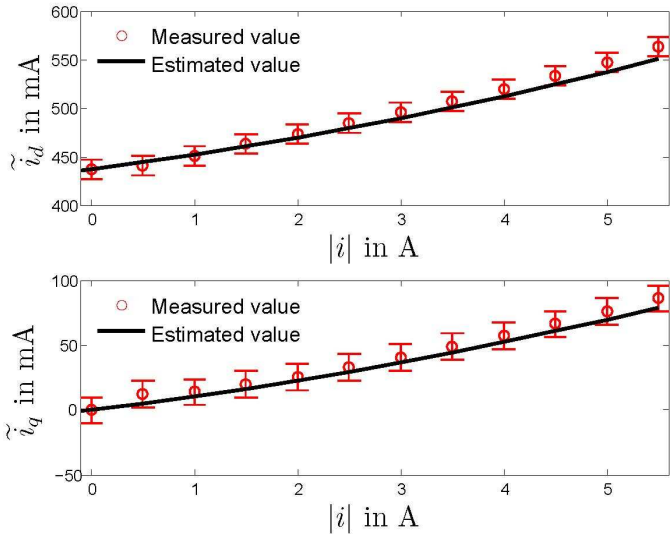}
   \label{fig:Id_i60_ud1}
 }
\caption{Measured values (circles) compared to model-predicted values (solid line) for a $60^\circ$ current angle.
   }\label{fig:i60_ud}
\end{figure}

\subsection{Validation}
The estimation procedure relies on~\eqref{eq:saturation1}--\eqref{eq:saturation4q}, with either $\bar i_d\neq$ or , i.e. current vectors with angles $0^\circ,90^\circ,180^\circ,270^\circ$. To check the validity of the model tests were conducted with current vectors with various angles and magnitudes on the whole operating ($|i|=\sqrt{i_d^2+i_q^2}$ ranging from $0A$ to $2A$ with a $0.3A$ increment for the IPM, and from $0A$ to $5.5A$ with a $0.5A$ increment for the SPM).
Fig.~\ref{fig:i60_ud} shows for instance the results for a $60^\circ$ current angle; there is a good agreement between the measured values and those predicted by the model.

\begin{figure}[t!]
\centerline{\includegraphics[width=\columnwidth]{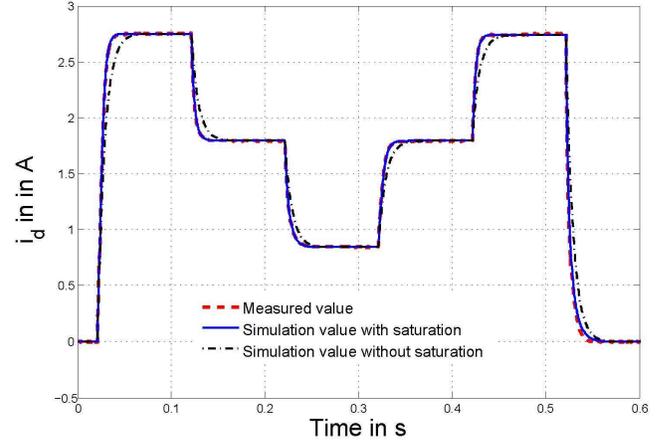}}
\caption{Time response of $i_d$ to large step voltages in $u_d$ (IPM).}\label{fig:step}
\end{figure}

\begin{figure}[t!]
\centerline{\includegraphics[width=\columnwidth]{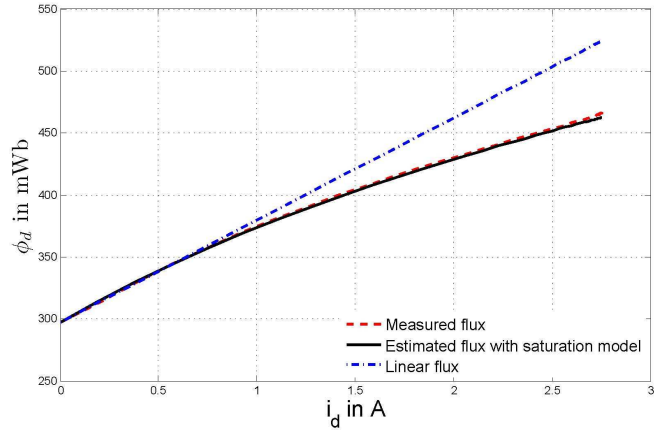}}
\caption{Saturation curve $\phi_d-i_d$ (IPM).}\label{fig:phididm}
\end{figure}

As a kind of cross-validation we also examined the currents time responses to large voltage steps. Fig.~\ref{fig:step} shows the good agreement between the measurements and the time response obtained by simulating the model with the estimated parameters; it also shows the differences with the simulated response when the saturation effects are omitted.
Fig.~\ref{fig:phididm} shows the good agreement also between the ``measured'' flux values (i.e. obtained by integrating the measured currents and voltages) and the flux values obtained by simulation.
%

\section{Conclusion}
 A simple parametric magnetic saturation model for the PMSM with a simple identification procedure based on high-frequency voltage injection have been introduced. Experimental tests on two kinds of PMSM (IPM and SPM) demonstrate the relevance of the approach. This model can be fruitfully used to design a sensorless control scheme at low velocity.

\end{document}